\documentclass[10pt,amscd]{amsart}
\usepackage{latexsym}
\usepackage{amsmath}
\usepackage{amssymb}
\usepackage[all]{xy}

\newcommand{\C}{\mathbb{C}}
\newcommand{\R}{\mathbb{R}}
\newcommand{\F}{\mathbb{F}}

\newcommand{\Z}{\mathbb{Z}}
\newcommand{\Q}{\mathbb{Q}}

\newcommand{\image}{{\rm Im}}

\newtheorem{thm}{Theorem }[section]
\newtheorem{prop}[thm]{Proposition}
\newtheorem{lem}[thm]{Lemma}

\newtheorem{cor}[thm]{Corollary}
\newtheorem{exmp}[thm]{Example}
\newtheorem{rem}[thm]{Remark}

\begin{document}
\title{A rational obstruction to be   a  Gottlieb map} 

\author{Toshihiro YAMAGUCHI}
\date{}
\address{Faculty of Education, Kochi University, 2-5-1, Kochi,780-8520, JAPAN}\email{tyamag@kochi-u.ac.jp}

\maketitle

\begin{abstract}
We investigate  {\it Gottlieb map}s,
which are  maps $f:E\to B$
that induce the maps between 
the  Gottlieb  groups
$\pi_n (f)|_{G_n(E)}:G_n(E)\to G_n(B)$ for all $n$,
from a rational homotopy theory point of view.
%Put $f_{\Q}$ to be  the rationalization of  $f$.
%We show that a map $f_{\Q}$ is a Gottlieb map
% a map $f:E\to B$  
 %$\pi_n (f_{\Q}):G_n(E_{\Q})\to G_n(B_{\Q})$ for all $n$
 %if $H^*(j;\Q):H^*(E;\Q)\to H^*(X;\Q)$ is surjective
 %for the homotopy fibration $X\overset{j}{\to} E\overset{f}{\to} B$.
We will define the   obstruction group $O(f)$
to be a Gottlieb map and a numerical invariant $o(f)$.
It naturally deduces   a  relative  splitting 
of  $E$ in certain cases. %in  the homotopy fibration$X\to E\to B$.
%which deduces a necessary condition to be a Gottlieb map.
We  also illustrate several  rational examples of Gottlieb maps and non-Gottlieb maps
by using derivation arguments in  Sullivan models.
%Our tool is the derivations of a Sullivan's model.
\end{abstract}
\section{Introduction}

The $n$th Gottlieb group (evaluation subgroup of homotopy group)
$G_n(B)$ of a path connected
CW complex $B$ with basepoint $*$ is the subgroup of the 
$n$th homotopy group $\pi_n(B)$ of $B$
consisting of homotopy classes of based maps $a:S^n\to B$ such that
the wedge $(a|id_B):S^n\vee B\to B$ extends
to a map $F_a:S^n\times B\to B$ \cite{G}.
%Note $G_n(Y)\subset G_n(Y,X;f)$ in general and $G_n(X)=G_n(X,X;id_X)$.
 The Gottlieb group is a very interesting homotopy invariant 
 (e.g., see \cite{O})
 but the calculations are difficult  even for spheres \cite{GM}.
 It is well known that 
 the Gottlieb group fails 
 to be a functor since, generally, a based map $f:E\to B$ does not yield a homomorphism
$\pi_*(f)|_{G_*(E)}:G_*(E)\to G_*(B)$ for $\pi_* (f):\pi_*(E)\to \pi_*(B)$.
For example,
$i: E=S^1\hookrightarrow S^1\vee S^1=B$
does not induce $\pi_1(i)|_{G_1(E)}:G_1(E)\to G_1(B)$
since $G_1(S^1)=\Z$ \cite[Theorem 5.4]{G} 
but $G_1(S^1\vee S^1)=0$ \cite[Theorem 3.1]{G}.
Recall that a space
$B$ is said  to be a {\it Gottlieb space }(or simply {\it G-space} in this paper) if
 $G_n(B)=\pi_n(B)$ for all $n$. 
 For example, an H-space is a G-space.
 It is  interesting to consider when is a space
 a G-space
 \cite{G},\cite{Si},\cite{La}.
In this paper, we will give a similar definition for a  map
and consider when is a map  such a map.
\vspace{2mm}

\noindent
{\bf Definition A.}
 If a map $f:E\to B$ induces $\pi_n (f)G_n(E)\subset G_n(B)$
 for all %$\pi_n (f):\pi_n(E)\to \pi_n(B)$  
 $n$,
  we call it   a {\it Gottlieb map} (or simply {\it G-map} in this paper). 
  \vspace{2mm}

We note some sufficient conditions to be a  G-map.
If $B$ is a G-space, any map $f:E\to B$ is a G-map. 
So `G-map' is a natural  generalization of `G-space'.
%and   is different in nuance with `H-map'  for `H-space'.
When
$E=S^n$, a G-map $f$  is an  $n$th Gottlieb element of $B$; i.e.,
$[f]\in G_n(B)$.
Also the projection $S^{d(n+1)-1}\to \F P^n$ for $d=\dim_{\R}{\F}$
 is a G-map under a certain condition \cite{GM2}.
 Here $\F P^n$ is the $n$-projective space over 
 $\F = \R , \C , {\mathbb{H}}$.
 If a map $f$ is  homotopic to the constant map; i.e., 
 $f\simeq *$,
 then it is a G-map.
 Put
$X\overset{j}{\to} E\overset{f}{\to} B$ the homotopy fibration 
where $X$ is the homotopy fiber of $f$.
Note that $f$ is a G-map if the fibration is fibre-homotopically trivial.
Also the connecting map
$\partial :\Omega B\to X$
is a G-map \cite{G2}.

The definition of $G_n(B)$
is generalized by replacing the identity
by an arbitrary based map   $f:E\to B$ \cite{WK}.
The $n$th evaluation subgroup $G_n(B,E;f)$ of the map $f$
is the subgroup of $\pi_n(B)$ 
 for the evaluation map 
$map(E,B;f)\to B$.
 It is 
represented by  maps $a:S^n\to B$ such that
 $(a|f):S^n\vee E\to B$ extends
to a map $F_a:S^n\times E\to B$. Put $G(Y)=\oplus_{i>0}G_i(Y)$ for
a space $Y$ 
and $G(B,E;f)=\oplus_{i>0}G_i(B,E;f)$.
From the definitions,  there is a map
 $\pi_*(f):G(E)\to G(B,E;f)$ %restricted from $\pi_n(f)$,
  and $G(B,E;f)\supset G(B)$.
  Therefore, the following is obvious.
  \begin{lem}
If $G(B,E;f)\subset G(B)$, %if $G(B,E;f)=G(B)$,
then $f:E\to B$ is a G-map.
\end{lem}
So if $f$ has a right homotopy inverse, 
$f$ is a G-map 
\cite[Proposition 1-4]{G}(\cite[Remark 3]{Y}).
For example, since
the free loop fibration  $\Omega X\to LX\overset{f}{\to} X$
has a section, the evaluation map $f$ is a G-map.
See \S 3 for the other sufficient conditions. 
%If $G(B)=G(B,E;f)$,$f$ is G-map. But 
%Even if $G(B) \neq G(B,E;f)$,$f$ may be a G-map (see Example 2 in \S 4). 
   
%In this paper, it will be useful
%to consider the (rational) homotopy fibration of a map.
%See \S 3.
%For the path fibration $\xi :\Omega X\to PX\to X$,$\partial$ satisfies $(i)$ for any $B$-fibration $\eta :B\to B'\to X$by $f':=*$. 

Suppose $E$ and $B$ have the homotopy types  of   nilpotent 
CW complexes.
Put $e_B:B\to B_{\Q}$ and $f_{\Q}=e_B\circ f
:E\to B_{\Q}$  to be the rationalizations of $B$
and $f:E\to B$, respectively \cite{HMR}.
Then $\pi_n(B_{\Q})\cong \pi_n(B)_{\Q}:=\pi_n(B)\otimes \Q$ for 
$n>1$.
By the universality of  rationalization, 
$f_{\Q}$ is equivalent to $\tilde{f_{\Q}}:
E_{\Q}\to B_{\Q}$, often we do not distinguish from $f_{\Q}$
in this paper.
%Put  $map (E,B;f)$ is the space of continuous (unbased) maps from $E$ to $B$ homotopic to $f$.
When $E$ is a finite complex, %${e_B}_*:map (E,B;f)\to map(E,B_{\Q};f_{\Q})$
 %is a rationalization \cite[Theorem II.3.1]{HMR}  and 
 $G_n(E_{\Q})\cong G_n(E)_{\Q}:=G_n(E)\otimes \Q$
 \cite{La2}
and $G(E_{\Q})$ is oddly graded \cite{FH}.
 Recall that $B_{\Q}$ is a G-space if and only if it is an H-space.
But it seems difficult to search a useful  necessary  and sufficient  condition
to be a  rationalized G-map. %(compare with that of a  cyclic map \cite{LS2}).
 If a map $f:E\to B$ induces $\pi_* (f_{\Q})G(E)\subset G(B_{\Q})$
 or $\pi_* ({f_{\Q}})G(E_{\Q})\subset G(B_{\Q})$,
  we call $f$ a {\it rational Gottlieb map} (or simply {\it r.G-map}).
  Of course, a G-map between nilpotent spaces is an r.G-map.
   For a map,
we can define an obstruction group:
\vspace{2mm}

%from the commutative diagram$$
%\xymatrix{   G_n(E)\ar[d]^{F_n}& G_n(B)\ar[l]_{\ \ \ \ \pi_n(f)}\ar[ld]^{inc.}\\  G(B,E;f) &     $$	  

%\begin{thm} If the homotopy fibration is $\Q$-TNCZ,
%$f$ is a r.G-map.\end{thm}

   %since the rationalization preserve Gottlieb group.
%As the main theorem of this paper, we give a sufficient condition  for   an r.G-map. 

\noindent
{\bf Definition B.}
The {\it $n$th obstruction group of a map $f:E\to B$ to be a  G-map} 
is given by 
${\it O}_n(f):=\image (\ \overline{\pi_n(f)|_{G_n(E)}}\ )
\subset \pi_n(B)/G_n(B)$.
Namely,
$$O_n(f):=\image \ ( \ G_n(E)\overset{\pi_n (f)}{\to} \pi_n(B)
 {\twoheadrightarrow} \pi_n(B)/G_n(B)\ ).$$
 {Also put 
${\it O}(f):=\oplus_{n>0} {\it O}_n(f)$
and denote $\dim O(f_{\Q})$ as $o(f)$.}
\vspace{2mm}

 Roughly speaking,
 $O_n(f)$ is of ``non-Gottlieb'' elements 
  in $\pi_n(B)$ yet ``Gottlieb'' in $\pi_n(E)$. 
  Recall that
$G_1(B)$ is contained in the center of $\pi_1(B)$ \cite[Corollary 2.4]{G}.
We have 
%${\it O}(f)$ is an oddly graded subgroup of $\pi_*(B)/G(B)$ and
that ${\it O}(f)=0$ if and only if $f$ is a G-map.
%If $E$ is  simply connected,${\it O}(f_{\Q})$ is oddly graded 
%since $G(E_{\Q})$ is oddly graded \cite{FHT}. 
%and${\it O}(f_{\Q})=0$ if and only if $f$ is a r.G-map.\\
Note that $o(f)$ is a  numerical rational homotopy invariant of a map
with $0\leq o(f)\leq min\{ \dim G(E_{\Q}), \dim \pi_*(B)_{\Q}-\dim G(B_{\Q})\}$
and it is a measure of the rational non-triviality of the homotopy fibration $X\to E\to B$.
If $f:Y\to Z$ and $g:X\to Y$   
are G-maps, then the composition $f\circ g:X\to Z$ is a G-map 
from $\pi_*(f)\circ \pi_*(g)=\pi_*(f\circ g)$.
It induces that $o(f\circ g)=0$ if 
$o(f)=0$ and $o(g)=0$.
It is generalized as  

%\begin{thm}Suppose that $E$ and $B$ are simply connected  CW complexes with rational homology of finite type.Then there is a map $\Psi :G(B_{\Q},E_{\Q};f_{\Q})\to \Ker H^*(f;\Q)$with $\Ker \Psi=G(B_{\Q})$.Therefore$O(f_{\Q})\subset \Ker H^*(f;\Q)$.\end{thm}

\begin{thm}For any maps $f:Y\to Z$ and $g:X\to Y$ between simply connected
  complexes of finite type,
 there is an inequality: 
$o(f\circ g)\leq o(f)+o(g).$
\end{thm}

%Recall Oprea's splitting theorem of the fiber of a fibration.
%Let $F'\to E'\to B'$ be a fibration.
%Then there is a subproduct $K\subset \Omega B'$ such that
%$F'\simeq {\mathcal F}\times K$ and $H^*(K)\cong Im (\partial^*:H(F')\to H^*(\Omega B')$.

%Recall for the LS-category of a map (\cite[page 35]{c}),
 %there is  an inequality:
%$cat(f\circ g)\leq cat(f)\cdot cat(g)$,
%which is induced from the definition.
%Note that $cat(id_X)=cat(X)$,
%which is the LS-category of a space $X$ (\cite{c}).
%But $o(id_X)=0$ for any $X$.
%Also $cat(c)=o(c)=0$ for a constant map $c$.

Notice that 
an element of $O(f_{\Q})$ is represented  by a map from the  product of  rationalized spheres $
a: K=S^{n_{1}}_{\Q}\times \cdots
\times S^{n_{k}}_{\Q}\to B_{\Q}$ by certain compositions
(\cite[p.494]{FL}).
%Suppose that 
%$a_i: S^{n_{i}}_{\Q}
%\to B_{\Q}$   ($1\leq i\leq k$)  are odd-spherical generators of $O(f_{\Q})$.
Suppose that $a$ makes odd-spherical generators.
Then there are a rational space $B_a$ and  a fibration 
$B_a\overset{j_a}{\to} B_{\Q}\overset{p_a}{\to} K$
given by the KS-extension
$(\Lambda (w_1,..,w_k),0)\to (\Lambda W,d_B)\to (\Lambda W_k,\overline{d_B})=M(B_a)$ 
with $|w_i|=n_i$,
$w_i^*=a|_{S_{\Q}^{n_i}}$ and $W=\Q\{w_1,..,w_k\}\oplus W_k$ in Sullivan's model theory \cite{FHT}
(see \S 2),
which satisfies $p_a\circ a\simeq id_K$.
Put $X_{\Q}\overset{g}{\to} F\overset{f_a}{\to} B_a$ the pull-back fibration of
$X_{\Q}\to E_{\Q}\overset{f_{\Q}}{\to} B_{\Q}$
by  $j_a$.
We  call it  the {\it pull-back fibration associated  to} $a$.
%Put $g:X_{\Q}\to F$ the fibre-inclusion of $f_k$.
Oprea's  homotopical
splittings   of rational spaces (\cite{O2}, \cite{O3}, \cite{H2},
\cite{FL})
implies
the following result. % above  obstruction seems  closely related to them.

\begin{thm}Let  $f:E\to B$ be a map between simply connected
  complexes of finite type 
%where $\dim H^*(E;\Q )<\infty$
  and $X\overset{j}{\to} E\overset{f}{\to} B$ the homotopy fibration.
If a map  $a : K=(S^{n_{1}}\times \cdots
\times S^{n_{k}})_{\Q}\to B_{\Q}$ of  odd-spherical generators of $B_{\Q}$
 represents an element  in
  $O(f_{\Q})$,
 then the  fibre-inclusion $g: X_{\Q}\to F$ of  the pull-back fibration associated  to $a$
induces a  splitting  $\psi_{f,a} : E_{\Q}\simeq F\times K$ such that 
 $$   
\xymatrix{ X_{\Q}\ar[d]_{g}\ar[r]^{j_{\Q}}&
 E_{\Q}\ar[d]_{\simeq}^{\psi_{f,a}}\\
F\ar[r]_{i_1}& 
 F\times K
\\
}\ \ and\ \ 
 \xymatrix{E_{\Q}\ar[r]^{f_{\Q}}&
 B_{\Q}\\
F\times K\ar[u]_{\simeq}^{\psi_{f,a}^{-1}}& 
 K\ar[u]_a\ar[l]^{\ \ \ i_2}
\\
}$$
with $i_1(x)=(x,*)$ and $i_2(x)=(*,x)$
homotopically commute.
Moreover, this splitting does not come  from that of $B_{\Q}$; i.e.,
 the maps $a_i:S^{n_i}_{\Q}\hookrightarrow K\overset{a}{\to}  B_{\Q}$
 {\it cannot} be extended to $S^{n_i}_{\Q}\times B'_i\simeq  B_{\Q}$
for any  space $B'_i$.

Conversely,  if there exists such a splitting 
$\psi_{f,a} : E_{\Q}\simeq F\times K$ for a map $f:E\to B$,
then the map $a:K\to B_{\Q}$ of odd-spherical generators represents
an element of $O(f_{\Q})$,
in particular 
 $k\leq o(f)$.
\end{thm}

Thus, if a map $f_{\Q}:E_{\Q}\to B_{\Q}$ is a G-map, then there exists no above splitting
$\psi_{f,a}$
of $E_{\Q}$.
That is a necessary condition to be a G-map but is  not sufficient (see Example 3.2(4)). 
Notice that  Oprea \cite[Theorem 1]{O2}, \cite[(RFDT)]{O3} gives
 a rational decomposition of  the 
fibre $X$ of a fibration $X\to E\to B$ (see Remark \ref{O}).
Also Halperin \cite[Lemma 1.1]{H2}
and F\'{e}lix-Lupton \cite[Theorem 1.6]{FL} 
(when we restrict their 
 generalized evaluation map \cite[Definition 1.1]{FL}  to 
$a:K\to E_{\Q}$ itself) give a rational decomposition of  a space $E$
and
 our theorem seems a  relative one of it.

Though Definition A is  defined for all  connected based CW complexes, 
we focus on simply connected  CW complexes $E$ with rational homology of 
finite type %(with the exception of the circle) 
 with $\dim G(E_{\Q})<\infty$
when we consider  rational homotopy types (Sullivan minimal models).
 We do not distinguish between a map and the homotopy class that it represents.
Our tool is the derivations 
(\cite{FH}, \cite{LS4},\cite{LS},\cite{Su}) of Sullivan models \cite{Su},
which are prepared in \S 2.
So we assume that the reader is familiar with the basics of rational homotopy
theory \cite{FHT}.
 We   see a property  of $O(f_{\Q})$ in Lemma 2.3 
 and prove Theorem 1.2 and Theorem 1.3 in \S 2.
We  will illustrate some  rational examples in \S 3,
in which we note  examples of r.G-maps
which do not satisfy Lemma 1.1 %$G(B_{\Q},E_{\Q};f_{\Q})=G(B_{\Q})$
 in Example 3.3.
Also we mention   interactions  with Gottlieb trivialities \cite{LS} in Remark 3.4,
 cyclic maps \cite{V}(\cite{LS2}) in Example 3.6 and {\it W-map}s (see Definition C) in Example 3.7.

%They may indicate the difficulty to characterize the class of G-maps.
%In this paper,  
%Even for the rational  cases,  our  G-maps may seem   vague.

%We note that there is an approach of function space to a non-Gottlieb elements, % detected by Whitehead products,without using derivation arguments \cite{HKO}. 

%Remark that various suitations on our G-map can occur even in this restriction.

\section{Derivations of Sullivan models}
We use the {Sullivan minimal model} $M(Y)$ 
 of a nilpotent space $Y$ of finite type.
It is a free $\Q$-commutative differential graded algebra (DGA) 
 $(\Lambda{V},d)$
 with a $\Q$-graded vector space $V=\bigoplus_{i\geq 1}V^i$
 where $\dim V^i<\infty$ and a decomposable differential;
 i.e., $d(V^i) \subset (\Lambda^+{V} \cdot \Lambda^+{V})^{i+1}$
 and $d \circ d=0$.
Here  $\Lambda^+{V}$ is 
 the ideal of $\Lambda{V}$ generated by elements of positive degree.
Denote the degree of a homogeneous  element $x$ of a graded algebra as $|{x}|$
and the $\Q$-vector space of basis $\{v_i\}_i$ as $\Q\{v_i\}_i$. 
Then  $xy=(-1)^{|{x}||{y}|}yx$ and $d(xy)=d(x)y+(-1)^{|{x}|}xd(y)$. 
A map $f:X\to Y$ has a minimal model
which is a  DGA-map $M(f):M(Y)\to M(X)$.
Notice that $M(Y)$ determines the rational homotopy type of $Y$.
Especially there is an isomorphism $Hom_i (V,\Q )\cong \pi_i(X)_{\Q}$.
See \cite{FHT} for a general introduction and the standard notations.

Let $A$ be a DGA $A=(A^*,d_A)$ with $A^*=\oplus_{i\geq 0}A^i,\ A^0=\Q$, $A^1=0$  and 
the augmentation $\epsilon:A\to \Q$.
Define  
 $Der_i A$  the vector space of self-derivations of $A$
decreasing the degree by $i>0$,
where  $\theta(xy)=\theta(x)y+(-1)^{i|x|}x\theta(y)$
for $\theta\in Der_iA$. 
We denote  $\oplus_{i>0} Der_iA$ by
$DerA$.
The boundary operator $\delta : Der_* A\to Der_{*-1} A$
is 
defined by $\delta (\sigma)=d_A\circ \sigma-(-1)^{|\sigma |}\sigma\circ d_A$.
%Then $(Der A,\delta )$ is a differential graded Lie algebra (DGL)
%by $[\theta, \theta']=
%\theta\circ \theta'-(-)^{|\theta ||\theta'|}\theta'\circ \theta$.
For a DGA-map $\phi:A\to B$,
define a $\phi$-derivation of degree $n$ to be a linear map
$\theta :A^*\to B^{*-n}$
 with $\theta(xy)=\theta(x)\phi (y)+(-1)^{n|x|}\phi (x)\theta(y)$
 and 
$Der(A,B;\phi)$ the vector space of $\phi$-derivations. 
The boundary operator $\delta_{\phi}: Der_* (A,B;\phi )\to Der_{*-1} (A,B;\phi )$ is defined by $\delta_{\phi}
 (\sigma)=d_B\circ\sigma-(-1)^{|\sigma |}\sigma\circ d_A$. 
Note $Der_* (A,A;id_A)=Der_*(A)$.
For 
$\phi:A\to B$, %=(\Lambda W,d_B)$,
the composition with $\epsilon':B\to \Q$ 
induces a chain map $\epsilon'_*:
 Der_n(A,B;\phi )\to Der_n(A,\Q;\epsilon )$.
 For a minimal model $A=(\Lambda Z,d_A)$, define 
$G_n(A,B;\phi):=\image (H(\epsilon'_*):H_n(Der(A,B;\phi))\to Hom_n(Z,\Q)).$
%Define the $n$th rationalized Gottlieb group by
% 
%That is, if $\delta(v^*+\theta)=0$ for $\theta\in Der (\Lambda W\otimes \Lambda V)$ with $\theta(v)=0$.
Especially  $G_*(A,A;id_A)=G_*(A)$. 
Note that $z^*\in Hom (Z,\Q)$ ($z^*$ is the dual  of the basis element $z$) is in 
 $G_n(A,B;\phi )$
  if and only if 
  $z^*$ extends to a derivation $\theta\in Der(A,B;\phi )$
  with $\delta_{\phi} (\theta)=0$.
%  This paper  strongly depends on the following  result.
%For example, see \cite[p.392-393]{FHT}.

\begin{thm}\cite{FH},\cite{LS4},\cite{Su}\ \ When $E$ and $B$ are simply connected, 
$ G_n(B_{\Q},E_{\Q};f_{\Q})$
$\cong G_n(M(B),M(E);M(f))$, 
in particular 
$ G_n(B_{\Q})\cong G_n(M(B))$.
\end{thm}

 Let $\xi: X\overset{j}{\to} E\overset{f}{\to} B$ be a fibration.
Put $M(B)=(\Lambda W,d_B)$. 
  Then the model (not minimal in general) of $E\to B$ is given by 
a KS(Koszul-Sullivan)-extension $(\Lambda W,d_B)\to (\Lambda W\otimes \Lambda V,D)$
 with  $D|_{\Lambda W}=d_B$ and 
% Given a map $f:X\to Y$, let $map(X,Y;f)$ denote the path component of $f$. When $X$ and $Y$ are simply connected with $X$ a finite complex, it is known that $\pi_n(map(X,Y;f))\otimes \Q\cong H_n(Der (M(Y),M(X);M(f)))$  for $n\geq 2$ [LS].which gives 
a DGA-commutative diagram
{\small $$\xymatrix{(\Lambda W,d_B)\ar@{=}[d] \ar[r]^{i}& (\Lambda W\otimes \Lambda V,D)\ar[d]^\psi_{\simeq} \ar[r]^{} &  (\Lambda V,\overline{D})\ar@{=}[r]
\ar[d]^{\cong}&(\Lambda V,d)\\
M(B)\ar[r]^{M(f)} & M(E) \ar[r]^{M(j)}&  M(X),&\\
}$$}
\noindent
where   %$(\Lambda V,\overline{D})=M(X)$ for the homotopy fiber $X$ of $f$ .
     `$\simeq$'  means to be quasi-isomorphic \cite[{\S 15}]{FHT}.
   Then $G_n(M(B),M(E);M(f))= G_n((\Lambda W,d_B),((\Lambda W\otimes \Lambda V,D);i)$.
   In this paper, we consider the models of  r.G-maps mainly 
   in KS-extensions.

\begin{exmp}{\rm
In general,
$\psi$ is not a DGA-isomorphism.
For example,
put $M(E)=M(S^3)=(\Lambda (x),0)$
and $M(B)=(\Lambda (w_1,w_2,u),d_B)$ with $|w_i|=3$, 
$d_Bw_i=0$ and  
$d_Bu=w_1w_2$.
Suppose that
a map $f:S^3\to B$
satisfies 
$M(f)(w_1)=x$, $M(f)(w_2)=0$, $M(f)(u)=0$.
Then $(\Lambda V,\overline{D})=(\Lambda (v_1,v_2),0)$
 with $|v_1|=2$ and $|v_2|=4$
and $\psi :(\Lambda (w_1,w_2,u,v_1,v_2),D)\to (\Lambda x,0)$
is given by 
$Dv_1=w_2$, $Dv_2=u+w_1v_1$, $\psi (w_1)=x$ and the others to zero.
It is quasi-isomorphic but not isomorphic.
%}Then$X_{\Q}\simeq K(\Q, 2)\times K(\Q, 4)$.
}\end{exmp}

  %If it is not an r.G-map,there is an element $\alpha\in \Lambda W\otimes \Lambda V$ with $D\alpha\neq 0\in \Lambda W$ from Proposition 2.4.Then $\alpha=w\beta+\gamma$ where $w\in \Lambda W$ and $\beta$ is a non-exact $d$-cocycle of $\Lambda V$. Since $D\beta\neq 0$, $H^*(X;\Q)\ni [\beta ]\not\in In\ j^*$.QED\\

%Especially, if $G(B_{\Q})\neq G(B_{\Q},E_{\Q};f_{\Q})$,ere exist such non-zero elements $g_i$.  It will be often useful for considering whether or not $G(B_{\Q})= G(B_{\Q},E_{\Q};f_{\Q})$.

From Theorem 2.1 and Definition B for a map $f:E\to B$, we have

 \begin{lem} For  $W=\Q\{w_i\}_{i\in I}$ where $\pi_* (B)_{\Q}
 =Hom (W,\Q)$ with $|w_i|\leq |w_j|$ if $i<j$,
 put $I':=\{ i\in I|[w_i]\neq 0 \ \ in \ \ H^*(  
W\oplus  V,Q(D))\}$.
 Then
 there is an isomorphism
  $$O(f_{\Q})\cong \Q\{w_i^*,
\ {i\in I'} |\ \  w_i^* \mbox{ satisfies (i) and (ii) }\}$$
 where (i)
$ \delta_E(w_i^*+\sigma)=0$ for  some $\sigma\in Der(\Lambda 
W\otimes \Lambda V,\delta_E) $ with $\sigma (w_j)=0$ for  $j\leq i$\\
and 
(ii) $ \delta_B(w_i^*+\tau )\neq 0$ for  any $\tau\in Der(\Lambda W,\delta_B) $
with $\tau (w_j)=0$ for  $j\leq i$.
%$w^*\not\in G(B)_{\Q}$.
 \end{lem}

Here $Q(D)$ is the linear part of $D$.
 For example, $O(f)_{\Q}\cong \Q\{w_1^*\}$
  in  Example 2.2 since in particular $\delta_E ((w_1,1)-(u,v_1))=0$ (see Notation below). 
  
%  and $\delta_E ((w_2,1)+\sigma )\neq 0$ 
 % for any $\sigma\neq -(w_2,1)\in Der_{|w_2|}(\Lambda W\otimes \Lambda V)$.\\
   
   Theorem 1.2 follows from 
  \begin{prop}For any maps $f:Y\to Z$ and $g:X\to Y$ between simply connected
 spaces,
  there is an inclusion 
  $O(f_{\Q}\circ g_{\Q})\subset O(f_{\Q})\oplus O(g_{\Q})$. 
% Thus $\rank O(f\circ g)\leq \rank O(f)+\rank O(g)$. 
 \end{prop}

\noindent
{\it Proof.} \ Put a model  of $f\circ g:X\to Y\to Z$ as 
the  commutative diagram
%$$(\Lambda W,d_Z)\to (\Lambda W\otimes \Lambda V,D)\to (\Lambda W\otimes \Lambda V\otimes \Lambda U,D'),$$
$${\small \xymatrix{(\Lambda W,d_B)\ar@{=}[d] \ar[r]^{ }& (\Lambda W\otimes \Lambda V,D)\ar[d]_{\simeq} \ar[r]^{ } & 
 (\Lambda W\otimes \Lambda V\otimes \Lambda U,{D}')
\ar[d]^{\simeq}\\
M(Z)\ar[r]^{M(f)} & M(Y) \ar[r]^{M(g)}&  M(X),\\
}}$$
where $D|_{\Lambda W}=d_Z$
and $D'|_{\Lambda W\otimes \Lambda V}=D$. 
 For  $W=\Q\{w_i\}_{i\in I}$,
 $I'=\{ i\in I|[w_i]\neq 0 \ \ in \ \ H^*(  
W\oplus  V,Q(D))\}$ and 
$I'\supset I'':=\{ i\in I|[w_i]\neq 0 \ \ in \ \ H^*(  
W\oplus  V\oplus U,Q(D'))\}$,
from Lemma 2.3,
$$O(g_{\Q})\cap W^*= \Q\{w_i^*,
\ {i\in I''}   |\ \  w_i^* \mbox{ satisfies (i) and (ii) }\}$$
 where (i)
$ \delta_X(w_i^*+\sigma)=0$ for  some $\sigma\in Der(\Lambda 
W\otimes \Lambda V\otimes \Lambda U) $ with $\sigma (w_j)=0$ for $j\leq i$ 
and 
(ii) $ \delta_Y(w_i^*+\tau )\neq 0$ for  any $\tau\in Der(\Lambda W\otimes 
\Lambda V) $
with $\tau (w_j)=0$ for  $j\leq i$,
$$O(f_{\Q})=\Q\{w_i^*,
\ {i\in I'} 
  |\ \  w_i^* \mbox{ satisfies (iii) and (iv) }\}$$
 where (iii)
$ \delta_Y(w_i^*+\sigma)=0$ for  some $\sigma\in Der(\Lambda 
W\otimes \Lambda V) $ with $\sigma (w_j)=0$ for $j\leq i$ 
and 
(iv) $ \delta_Z(w_i^*+\tau )\neq 0$ for  any $\tau\in Der(\Lambda W) $
with $\tau (w_j)=0$ for  $j\leq i$, and 
$$O(f_{\Q}\circ g_{\Q})=\Q\{w_i^*,
\ {i\in I''} 
  |\ \  w_i^* \mbox{ satisfies (i) and (iv)  }\}.$$
Since
(ii) and (iii) contradict, 
we have $O(f_{\Q})\cap O(g_{\Q})=0$ in $W^*\oplus V^*$.
Also 
if $w_i^*\in O(f_{\Q}\circ g_{\Q})$ and $w_i^*\not\in O(f_{\Q})$,
then $w_i^*$ satisfies (i) but  not (iii).
Thus
$w_i^*\in O(g_{\Q})$.
\hfill\qed\\

\noindent
{\it Proof of Theorem 1.3.}\ \ Put the KS-extension of $f$ 
$(\Lambda W,d_B)\to (\Lambda W\otimes \Lambda V,D)$. 
For a sub-basis   $\{ w_{1},..,w_{k} \}$ of $W$,
 put $O(f_{\Q})\supset \Q\{w_{1}^*,..,w_{k}^*\}$ %($i_1< \cdots < i_k$)
with $|w_{i}|=n_{i}$ odd 
and $H^*(K;\Q)\cong \Lambda (w_{1},..,w_{k})$.
%Write simply $w_{i_1},..,w_{i_k}$ as $w_1,..,w_k$.
The assumption  induces $D(w_i)=d_B(w_i)=0$ for $i=1,..,k$.
From Lemma 2.3,
$ \delta_E(w_i^*)=\delta_E(\sigma_i)$ for  some $\sigma_i\in Der(\Lambda 
W\otimes \Lambda V) $.  % with $\sigma_i (w_j)=0$ for  $j\leq i$.
Put
$D_1=D$ and 
$D_{i+1}=\varphi_i^{-1}\circ D_i\circ \varphi_i
$ for $\varphi_i=id-\sigma_i\otimes w_i$ inductively for $i=1,..,k$,
which induce the changes of basis:
$$\varphi_i:
(\Lambda W_i\otimes \Lambda V,D_{i+1})\otimes (\Lambda (w_1,..,w_{i}),0)
\cong (\Lambda W_{i-1}\otimes \Lambda V,D_{i})\otimes (\Lambda (w_1,..,w_{i-1}),0)$$
for $W=W_i\oplus \Q\{w_1,..,w_{i}\}$ 
 \cite[Lemma 1.1]{H2} (the proof of \cite[Lemma A]{Y}).
Thus there is a DGA-isomorphism 
$$\varphi_1\circ \cdots \circ \varphi_k:
(\Lambda W_k\otimes \Lambda V,D_{k+1})\otimes (\Lambda (w_1,..,w_{k}),0)
\cong (\Lambda W\otimes \Lambda V,D).$$
The model of the pull-back 
$$\xymatrix{ B_a
\ar[d]_{j_a} &  F\ar[l]_{f_a}
\ar[d]\\
B_{\Q}& E_{\Q}\ar[l]^{f_{\Q}}
\\
}$$
is given by the push-out
$$\xymatrix{  (\Lambda W_k,\overline{d_B})
\ar[r]^{} &  (\Lambda W_k\otimes \Lambda V,\overline{D})\\
 (\Lambda W,d_B) \ar[r]\ar[u]& (\Lambda W\otimes \Lambda V,D)\ar[u]
\\
}$$
with $\overline{D}|_{\Lambda W_k}=\overline{d_B}$.
Notice that  $M(F) =(\Lambda W_k\otimes \Lambda V,\overline{D})\cong (\Lambda W_k\otimes \Lambda V,D_{k+1})$
and then the model of $g:X_{\Q}\to F$ is given by 
the projection $p:(\Lambda W_k\otimes \Lambda V,D_{k+1})
\to (\Lambda V,\overline{D_{k+1}})= (\Lambda V,d)$.
We have
the DGA-commutative diagrams
$$\xymatrix{  (\Lambda (w_1,..,w_{k}),0)\otimes(\Lambda W_k\otimes \Lambda V,D_{k+1})
\ar[d]_{\varphi_1\circ \cdots \circ \varphi_k}^{\cong} \ar[r]^{} &  (\Lambda W_k\otimes \Lambda V,D_{k+1})
\ar[d]^p\\
 (\Lambda W\otimes \Lambda V,D) \ar[r]& (\Lambda V,d)
\\
}$$
and
$$
\xymatrix{(\Lambda (w_1,..,w_{k}),0)&(\Lambda (w_1,..,w_{k}),0)\otimes(\Lambda W_k\otimes \Lambda V,D_{k+1})\ar[l]\\
(\Lambda W,d_B)\ar[u]^{M(a)}\ar[r]&(\Lambda W\otimes \Lambda V,D). \ar[u]_{(\varphi_1\circ \cdots \circ \varphi_k)^{-1}}^{\cong}\\
}$$
They are the models of the diagrams in Theorem 1.3.

The converse is given as follows.
The odd-spherical generators   $a_i:S^{n_i}_{\Q}\hookrightarrow K\overset{a}{\to}  B_{\Q}$
are not in $G(B_{\Q})$ from the assumption \cite[Lemma 1.1]{H2}.
On the other hand,
$\psi^{-1}_{f,a}|_{S^{n_i}_{\Q}}\in G(E_{\Q})$
from $\psi_{f,a}:E_{\Q}\simeq F\times S^{n_{1}}_{\Q}\times \cdots
\times S^{n_{k}}_{\Q}$.
Since $f_{\Q}\circ\psi^{-1}_{f,a}|_{S^{n_i}_{\Q}}\simeq a_i $, we have
 $a_i\in O(f_{\Q})$
from Definition B.
\hfill\qed\\

From Theorems 1.2, 1.3 and Proposition 2.4, we have
%Remark that the proof is similar to  arguments  in \cite{O2}, \cite{O3}, \cite{H2}  and   \cite{FL}.

\begin{cor}\label{A}
For maps $f:Y\to Z$ and $g:X\to Y$,
if there is a splitting 
$\psi_{f\circ g, a}:X_{\Q}\simeq F\times K$ as in Theorem 1.3,
where a map $a : K=(S^{n_{1}}\times \cdots\times S^{n_{k}})_{\Q}\to Z_{\Q}$ makes  odd-spherical generators of $Z_{\Q}$, 
then 
$k\leq o(f)+o(g)$. Also,
 % of $\dim H^*(X;\Q)<\infty$ and  $\dim H^*(Y;\Q)<\infty$,

\noindent
(i) Suppose
$O(f_{\Q}\circ g_{\Q})= O(f_{\Q})\oplus O(g_{\Q})$. 
If elements $a:K\to Y_{\Q}$ of $O(g_{\Q})$ and 
 $b:K'\to Z_{\Q}$ of $O(f_{\Q})$
make both odd-spherical generators,
then there is a decomposition 
$X_{\Q}\simeq F\times K\times K'$
for some rational space $F$.

\noindent
(ii) Suppose that $g$ is an r.G-map.  
If  there is a splitting $\psi_{f\circ g,a}:X_{\Q}\simeq F\times K$ as in Theorem 1.3,
then it deduces  a splitting $\psi_{f,a}:Y_{\Q}\simeq F'\times K$
for some rational space $F'$.
\end{cor}

\begin{rem}{\rm 
(1) %In theorem 1.3, we need an element $a$   to be spherical. 
%there may be a splitting
%even if $o(f)=0$.
Put $B$ the homogeneous space 
$SU(6)/SU(3)\times SU(3)$ ($SU(n)$ is a special unitary group), 
whose model is given by
$(\Lambda (x,y, v_1,v_2,v_3),d_B)$ with
$|x|=4$, $|y|=6$, $|v_1|=7$, $|v_2|=9$, $|v_3|=11$, $d_Bx=d_By=0$,
$d_Bv_1=x^2$, $d_Bv_2=xy$ and  $d_Bv_3=y^2$ \cite[p.486]{GHV}.
For a map $f:E\to B$ of the KS-extension 
$(\Lambda (x,y, v_1,v_2,v_3),d_B)\to (\Lambda (x,y, v_1,v_2,v_3,v),D)$
with $|v|=3$ and 
$Dv=x$, we have 
$o(f)=0$ but there is a splitting  $\psi_{f,a}: E_{\Q}\simeq F\times (S^7\times S^{9})_{\Q}$
for a map of  (non-spherical)  Gottlieb elements 
$a:(S^7\times S^{9})_{\Q}\to B_{\Q}$ and $F=S^6_{\Q}$.
We note 
$ (S^7\times S^{9})_{\Q}=K_4$
in Theorem \ref{D}  below.

(2) For a map $f:E\to B$,
if an element $a : K=(S^{n_{1}}\times \cdots
\times S^{n_{k}})_{\Q}\to B_{\Q}$
 of  $O(f_{\Q})$
 makes  odd-spherical generators of $B_{\Q}$, then 
 we see from the second diagram in the proof of Theorem 1.3 that
the pull-back fibration $X_{\Q}\to E'\to K$
of the homotopy fibration $X_{\Q}\to E_{\Q}\to B_{\Q}$ by $a:K\to B_{\Q}$
is fibre-homotopically trivial.
Indeed, the model is given by the push-out
 $$\xymatrix{  (\Lambda (w_1,..,w_{k}),0)
\ar[r]^{} &  (\Lambda (w_1,..,w_{k})\otimes \Lambda V,\overline{D_{k+1}})
\ar@{=}[r]
& (\Lambda (w_1,..,w_{k}),0)\otimes (\Lambda V,d)\\
 (\Lambda W,d_B) \ar[r]\ar[u]^{M(a)}& (\Lambda W\otimes \Lambda V,D).\ar[u]_{\overline{(\varphi_1\circ \cdots \circ \varphi_k)^{-1}}}&
\\
}$$

(3)  For a map $f:E\to B$,
suppose that $f_a:F\to B_a$  is the pull-back fibration associated  to
a map $a : K=(S^{n_{1}}\times \cdots
\times S^{n_{k}})_{\Q}\to B_{\Q}$ of odd-spherical generators of $O(f_{\Q})$.
Then $o(f_a)\leq o(f)-k$.
}
\end{rem}

For a fibration $\xi :X\overset{j}{\to} E\overset{f}{\to} B$ of rational spaces,
 there is a decomposition $G_n(E)= S_n\oplus T_n\oplus U_n
\subset G_n(X)\oplus GH_n(\xi )\oplus G_n(B,E;f)$ where
 $U_n:=\pi_n(f)(G_n(E))\subset G_n(B,E;f)$ 
\cite[Theorem A]{Y} and then
$O_n(f)= U_n/(G_n(B)\cap U_n)$.
Here  $GH_n(\xi ):=Ker (\ \pi_n(f):G_n(E,X;j)\to \pi_n(B)\ )\ /\ Im(\ \pi_n(j):G_n(X)\to G_n(E,X;j)\ )$
 is called as the $n$th Gottlieb homology group  of $\xi$ \cite{LW}, \cite{LS}.
From the manner of \cite[Theorem A]{Y},
we have 

\begin{thm}\label{D} For a fibration $\xi :X\overset{j}{\to} E\overset{f}{\to} B$ of rational spaces,
suppose that there is a decomposition $E\simeq F\times S$ where
a map $a:S=(S^{n_{1}}\times \cdots
\times S^{n_{k}})_{\Q}\to E$ makes  odd-spherical generators.
Then $S$ is uniquely decomposed as $S=
K_1\times K_2\times K_3\times K_4$
where $a|_{K_1}$ makes generators of $\pi_*(j)G(X)$,
$a|_{K_2}$ makes generators of $GH(\xi)$,
$f\circ a|_{K_3}$ makes generators of $O(f)$
and $f\circ a|_{K_4}$ makes generators of $G(B)$.
In particular,
$K_3=*$ if $f$ is a G-map and 
$K_2=K_3=*$ if $\xi$ is a trivial fibration.
\end{thm}

\begin{rem}\label{O}{\rm
 Recall  
Oprea's rational fibre decomposition theorem(\cite{O2},\cite{O3},\cite{O}):
{\it For a fibration $\xi :X\overset{j}{\to} E\overset{f}{\to} B$ of rational spaces
with finite betti  numbers,
there is a subproduct ${\mathcal K}\subset \Omega B$ and  a space ${\mathcal F}$ such that
$X\simeq {\mathcal F}\times {\mathcal K}$ and $H^*({\mathcal K})\cong Im (\partial^*:H^*(X)\to H^*(\Omega B))$}.
The space ${\mathcal K}$ is called the Samelson space of $\xi$.
If we apply this theorem to the  rationalized Hopf fibration $S^3_{\Q}\overset{j}{\to} S^7_{\Q}\to S^4_{\Q}$,
% $E_{\Q}\overset{f}{\to} B_{\Q}\to K(\Q, 6)$ induced from above,
  the Samelson space is the fibre $S^3_{\Q}$ itself.
But it can not be $K$ in 
Theroem 1.3 since $o(j)=0$
for the induced fibration $\Omega S^4_{\Q}\to S^3_{\Q}\overset{j}{\to} S^7_{\Q}$.
In general,  in Theorem \ref{D}  
for the induced fibration $\Omega B\overset{\partial}{\to} X\overset{j}{\to} E$,
we have
 $K_1\subset {\mathcal K}$
as a subproduct, $K_2=*$ and $K_i\cap {\mathcal K}=*$ for $i=3,4$.
}
\end{rem}

\noindent
{\bf Notation}(\cite[Definition 16]{Sa},\cite[p.314]{Su}).
 For a DGA-map $\phi :(\Lambda V,d)\to (\Lambda Z,d')$, 
 the symbol $(v,h)\in Der (\Lambda V,\Lambda Z;\phi )$
 means the $\phi$-derivation sending an element $v\in V$
to $h\in \Lambda Z$ and the other to zero.
Especially $(v,1)=v^*$.
The differential is given as 
$$\delta_{\phi}(v,h)\ =\
d'\circ(v,h)-(-1)^{|v|-|h|}(v,h)\circ d
\ =\  (v,d'h)-\sum_i\pm_i (u_i,\phi (\partial du_i/\partial v)\cdot h)
$$
for a basis $\{u_i\}$ of $V$.
If $\phi=M(f)$ or a KS-extension of $M(f)$, we denote $\delta_{\phi}$
simply  as $\delta_f$.
We  often use the symbol  $( *, *)$
in the following section. 

\section{ Examples }

Fix the KS-model of   a based map $f:E\to B$
  as a DGA-map
$i:(\Lambda W,d_B)\to (\Lambda W\otimes \Lambda V,D),$
  where 
$D|_W=d_B$ and 
$(\Lambda V,\overline{D})=(\Lambda V,d)=M(X)$ for the homotopy fiber $X$ of $f$. 
 %   Denote $I(S)$ the ideal in a certain  algebra generated by a set $S$.
 % }\end{exmp}

% \begin{lem}If $D(\alpha )\in \Lambda W\otimes \Lambda^+V$ for any 
%non-exact $d$-cocycle $\alpha$of $\Lambda V$, %with $\alpha\not\in \Lambda W$,
%then $f$ is an r.G-map.\end{lem}
  
 \vspace{0.3cm}
 
\begin{exmp}{\rm
  %(2) If $X$ has the rational homotopy type of 
 %a simply connected compact homogeneous space $G/H$ with $\rank G=\rank H$,
 %then any map   $f:E\to B$ with fiber $X$ is an r.G-map from Theorem 1.3
 %since $H^*(E;\Q)\cong H^*(X;\Q)\otimes H^*(B;\Q)$as vector spaces  \cite{ST}.

Suppose $\dim H^*(E;\Q)<\infty$.
If  $B$ is  pure; i.e.,
$\dim W<\infty$, $d_BW^{odd}\subset \Lambda W^{even}$ and $d_BW^{even}=0$,
 then any map $f:E\to B$ is an r.G-map.
In fact, 
since  $G(E_{\Q})$ has generators of  odd degrees \cite[Theorem III]{FH},
we have $\pi (f)_{\Q}:G(E_{\Q})=G_{odd}(E_{\Q})\to \pi_{odd}(B_{\Q})=G(B_{\Q})$.
In particular, a map whose target is a homogeneous space is an r.G-map.     }
\end{exmp}

%\vspace{0.2cm}

\begin{exmp}{\rm We note some   rational splittings  obtained from    non-r.G-maps.

(1)
 Put an odd spherical fibration $S^m\to E\overset{f}{\to} B$ where
$M(S^m)=(\Lambda (v),0)$
and $M(B)=(\Lambda (w_1,w_2,\cdots ,w_{2n},u),d_B)$  $(n>1)$
with $m=|v|=|w_1|+|w_2|-1$, $|w_1|,\cdots , |w_{2n}|$, $|u|$ odd. 
When $d_Bu=w_1w_2\cdots w_{2n}$ and $Dv=w_1w_2$,
we have $\delta (w_i,1)(u)=D(w_i,1)(u)+(w_i,1)D(u)=(w_i,1)Du=(w_i,1)(w_1\cdots w_{2n})=
(-1)^{i-1}w_1w_2\cdots \overset{\vee}{w_i} \cdots w_{2n}$.
Then
$\delta_E ((w_i,1)+ (-1)^{i}(u,vw_3\cdots \overset{\vee}{w_i} \cdots w_{2n}))=0$
for $i=3,\cdots ,2n$.
Thus ${\it O}(f_{\Q})\cong \Q\{w^*_3,\cdots , w^*_{2n}\}$
from Lemma 2.3;
i.e., $f$ is not an r.G-map.
There is 
a decomposition
$$E_{\Q}\simeq F\times K=F\times S_{\Q}^{|w_{3}|}\times \cdots \times 
S_{\Q}^{|w_{2n}|}$$
where
$F=F'\times S_{\Q}^{|u|}$ with
  $M(F')\cong (\Lambda (w_1,w_2,v),d')$
with $d'w_i=0$ and $d'v=w_1w_2$.

(2) %Even if $H^*(f):H^*(B;\Q)\to H^*(E;\Q)$ is injective,$f$ is not an r.G-map in general. 
Put
$M(X)=(\Lambda (v,v'),0)$
and $M(B)=(\Lambda (w_1,w_2,w_3,w_4,u),d_B)$
with $d_Bw_i=0$ and $d_Bu=w_1w_2+w_3w_4$ where $|w_i|$ are odd
($|v|$ and $|v'|$ are even).
Suppose that the differential of the model of
a map $f:E\to B$ with homotopy fibre $X$ is given by 
$D(v)=0$, $D(v')=w_2+w_3v.$
 Then we have $\delta_E((w_1,1)-(u,v')-(w_4,v))=0.$
 Thus ${\it O}(f_{\Q})\cong \Q\{w^*_1\}$ from Lemma 2.3;
i.e., $f$ is not an r.G-map.
There is 
a decomposition
$$E_{\Q}\simeq F\times K=F\times S_{\Q}^{|w_{1}|}$$
where $F=F'\times K({\Q},{|v|})$
with  $M(F')\cong (\Lambda (w_3,w_4,u),d')$
with $d'w_i=0$ and $d'u=w_3w_4$.

(3) Put $E=S^3$ and $M(B)=(\Lambda (w_1,w_2,u),d_B)$ with the map of Example 2.2.
Then $M(X)\cong (\Lambda (v_1,v_2),0)$ for 
the homotopy  fibre $X$ and we have $$E_{\Q}\simeq K=S^{|w_1|}_{\Q}\ \ \mbox{ and }\ \ F=*$$
in  Theorem 1.3.
Here $M(F)=
(\Lambda W_1\otimes V,D_2)=(\Lambda (w_2,u,v_1,v_2),D_2)$ with $D_2w_2=D_2u=0$, $D_2v_1=w_2$
and $D_2v_2=u$ (see the proof of Theorem 1.3).

(4)
Put $E=SU(6)/SU(3)\times SU(3)$, 
where
$M(E)=(\Lambda (x,y, v_1,v_2,v_3),d_E)$ with
$|x|=4$, $|y|=6$, 
$d_Ex=d_Ey=0$,
$d_Ev_1=x^2$, $d_Ev_2=xy$ and  $d_Ev_3=y^2$.
Put $M(B)=(\Lambda (w,x,y, v_1,v_2,v_3),d_B)$ with $|w|=3$,
$d_Bw=d_Bx=0$, 
$d_By=wx$, $d_Bv_1=x^2$, $d_Bv_2=xy+wv_1$, $d_Bv_3=y^2+2wv_2$
and the KS-extension of $f:E\to B$
is given by $(\Lambda (w,x,y, v_1,v_2,v_3),d_B)\to (\Lambda (w,x,y, v_1,v_2,v_3,v),D)$
with $|v|=2$, $Dv=w$ and $D=d_B$ for the other elements.
Then $O(f_{\Q})=\Q\{v_1^*,v_2^*\}$.
But $E_{\Q}$ can not non-trivially decompose; i.e.,
$K\simeq *$ if $E_{\Q}\simeq F\times K$,
%for some product of rationalized odd-spheres $K$,
from the DGA-structure of $M(E)$.
Thus 
 the splitting of Theorem 1.3 does not fold  for non-spherical generators of $B$ in general.
}
\end{exmp}

%\vspace{0.3cm}

\begin{exmp}{\rm
A map $f:E\to B$ may be an r.G-map even if 
$G(B_{\Q})\neq G(B_{\Q},E_{\Q};f_{\Q})$ (see Lemma 1.1). 

(1) The Hopf
map $f:S^3{\to} S^2$
 is a G-map and $G_n(S^2,S^3;f )=\pi_n(S^2)$ for all $n$ \cite[Example 2.7]{LS}
but $\pi_2(S^2)={\mathbb{Z}}\neq 0=G_2(S^2)$.

(2) Consider the pull-back fibration of the Hopf fibration $S^3\to S^7\to S^4$,
$S^3\to E \overset{g}{\to} B={\mathbb{C}}P^2$,
induced by the  map
$\mathbb{C}P^2\twoheadrightarrow S^4$
obtained by pinching out the 2-cell. 
Put $M(S^3)=(\Lambda v,0)$ and $M({\mathbb{C}}P^2)=(\Lambda (w,u),d_B)$
with $|w|=2$, $|u|=5$, $d_Bw=0$ and $d_Bu=w^3$.
Then
the KS-extension is given by
$(\Lambda (w,u),d_B)\to (\Lambda (w,u,v),D)\to (\Lambda v,0)$
with $Dv=w^2$.
Then $(\Lambda (w,u,v),D)\cong (\Lambda (w,v),D)\otimes (\Lambda u,0)$;
i.e.,
$E_{\Q}\simeq 
(S^2\times S^5)_{\Q}$.
Then $g_{\Q}$ is a G-map.
In fact,
for $G(E)_{\Q}=G_3(E)_{\Q}\oplus G_5(E)_{\Q}$,
$\pi_3(g)_{\Q}=0$ and 
$\pi_5(g)_{\Q}:G_5(E)_{\Q}=\Q\{u^*\}\cong G_5(B)_{\Q}$.
In this case,
 $G(B)_{\Q}=\Q\{u^*\}\subset \Q\{w^*,u^*\}=G(B,E;g)_{\Q}$
 from $\delta_g((w,1)-(u,v))=0$.

(3)
Put $M(X)=(\Lambda (v),0)$ and $M(B)=(\Lambda (w_1,w_2,w_3,w_4,u),d_B)$
with $d_Bw_i=0$ and $d_Bu=w_1w_2w_3w_4$ where the degrees are odd.
If  $D(v)=w_1w_2+w_3w_4$
in a KS-extension,
it is  an r.G-map by direct calculation.
For example,  $\delta_E(w_1,1)=(u,w_2w_3w_4)+(v,w_2)$
and $\delta_E((w_1,1)+\sigma )\neq 0$ for any derivation $\sigma
\neq -(w_1,1)$.
Thus $O(f_{\Q})=0$ from Lemma 2.3.
In this case,
 $G(B_{\Q})=\Q\{u^*\}$ but $G(B_{\Q},E_{\Q};f_{\Q})=\Q\{ w^*_1,w^*_2,w^*_3,w^*_4,u^*\}=\pi_*(B)_{\Q}$.
In fact, for example,
we have $w_1^*\in G(B_{\Q},E_{\Q};f_{\Q})$ from $\delta_f((w_1,1)-(u,vw_2))=0$.

(4) Put
$M(X)=(\Lambda (v,v'),0)$
and $M(B)=(\Lambda (w_1,w_2,w_3,w_4,u),d_B)$
with $d_Bw_i=0$ and $d_Bu=w_1w_2w_3w_4$ where the degrees are odd.
If 
 $D(v)=w_1w_2+w_3v'${ and }$D(v')=w_3w_4$
in a KS-extension,
by direct calculation, 
%$$\delta_E ((w_1,1)-(u,v'w_2))=(v,w_2)$$$$\delta_E ((w_2,1)+(u,v'w_1))=-(v,w_1)$$$$\delta_E ((w_3,1)-(u,vw_4))=(v,v')+(v',w_4)$$$$\delta_E ((w_4,1)+(u,vw_3))=-(v',w_3).$$
 we see $O(f_{\Q})=0$ from Lemma 2.3.
In this case, $G(B_{\Q})=\Q\{u^*\}$ 
but $G(B_{\Q},E_{\Q};f_{\Q})=\pi_*(B)_{\Q}$.}
\end{exmp}

%\noindent{\bf Question 1.} {\it If a fibration $X\to E\overset{f}{\to} B$ c-splits, is $f$ a r.G-map ? }\\

\begin{rem} {\rm
Put $h:B\to Baut_1X$ the classifying map of a fibration of finite complexes
$\xi :X\overset{j}{\to} E\overset{f}{\to} B$.
If the rationalized 
Gottlieb sequence \cite{LW},\cite{LS} deduces the short exact sequence
$0{\to} G_n(X)_{\Q}
\overset{\pi_n(j)_{\Q}}{\to} G_n(E,X;j)_{\Q}\overset{\pi_n(f)_{\Q}}{\to} \pi_n(B)_{\Q}\to 0$
%i.e., $G_n(E,X;j)_{\Q}=G_n(X)_{\Q}\oplus \pi_n(B)_{\Q}$, 
 for all $n>1$,
the fibration $\xi$ is said to be rationally  Gottlieb-trivial
 \cite{LS}.
It is a notion of the relative triviality of fibration, too.
Recall that $f:E\to B$ is  rationally Gottlieb-trivial
if and only if  $\pi_* (h)_{\Q}=0$ \cite[Theorem 4.2]{LS}.
On the other hand, 
$\pi_* (h)_{\Q}$
cannot determine whether $f$ is an r.G-map or not.
For example, the Hopf
bundle $S^1\to S^3\overset{f}{\to} S^2$
 (1) and the fibration  (4) of Example 3.3 
are not rationally  Gottlieb-trivial
since $\pi_* (h)_{\Q}\neq 0$ from \cite[Theorem 3.2]{LS},
 %:\pi_4(S^4)_{\Q}\to \pi_4(BS^3)_{\Q}$,
but they are r.G-maps.
Also for the fibrations of Example 3.2 (1) and   of Example 3.3 (2), (3),
we see 
$\pi_* (h)_{\Q}=0$ from \cite[Theorem 3.2]{LS}.
From the definition of $GH(\xi )$,
we see  $K_2=*$ in Theorem \ref{D} if $\xi$ is  rationally Gottlieb-trivial.
}
\end{rem}

\begin{exmp}{\rm 
(1) Consider the  homotopy pull-back diagram of rational spaces:
$$
\xymatrix{\small E'\ar[d]_{f'}\ar[r]^{g'}& E\ar[d]^f\\
B'\ar[r]_{g} & B\\
}$$
where  $M(f):(\Lambda W,d_B)\to (\Lambda (W\oplus v),D)$,
 $M(g):(\Lambda W,d_B)\to (\Lambda (W\oplus v'),D')$
and the  homotopy groups are oddly graded.
Suppose $M(B)=(\Lambda (w_1,\cdots ,w_{2n},u),d_B)$
with $d_Bw_i=0$ and  $d_Bu=w_1\cdots w_{2n}$ $(n\geq 3)$.

Put
$Dv=w_1\cdots w_4$ and $D'v'=w_1w_2$.
 Then
$o(g)=o(f\circ g')=2n-2$, 
$o( f)=2n-4$ and $o(g')=2$.
Then $o(f\circ g')=o(f)+o(g')$, especially 
 $O(f\circ g')=O(f)\oplus O(g')$.
Thus there is a decomposition
$$E'\simeq F\times K\times K'=F\times S^{|w_3|}_{\Q}\times\cdots \times  
S^{|w_{2n}|}_{\Q}$$
as in Corollary \ref{A} (i).
Here $M(F)=(\Lambda (w_1,w_2,v'),D')\otimes (\Lambda (v,u),0)$, $K=S^{|w_3|}_{\Q}\times  
S^{|w_{4}|}_{\Q}$ for $g'$
and $K'=S^{|w_5|}_{\Q}\times\cdots \times  
S^{|w_{2n}|}_{\Q}$ for $f$.
Also from 
$o( f')=0$,
the above decomposition deduces
$$B'\simeq F'\times S^{|w_3|}_{\Q}\times\cdots \times  
S^{|w_{2n}|}_{\Q}$$
as in Corollary \ref{A} (ii).
Here  $M(F')=(\Lambda (w_1,w_2,v'),D')\otimes (\Lambda u,0)$.

Put $Dv=w_1\cdots w_4$ and $D'v'=w_5w_6$.
Then
$o(g)=2n-2$, $o( f)=2n-4$, 
$o(g\circ f')=o(f\circ g')=2n-6$
and 
$o( f')=o(g')=0$.
Then  there is a decomposition 
 $$E'\simeq F\times S^{|w_7|}_{\Q}\times\cdots \times  
S^{|w_{2n}|}_{\Q}$$
and it deduces $$B'\simeq F'\times S^{|w_7|}_{\Q}\times\cdots \times  
S^{|w_{2n}|}_{\Q}\ \ \mbox{and}\ \ E\simeq F''\times S^{|w_7|}_{\Q}\times\cdots \times  
S^{|w_{2n}|}_{\Q}$$as in Corollary \ref{A} (ii).
Here $M(F)=(\Lambda (w_1,\cdots ,w_6,v,v'),D'')\otimes (\Lambda u,0)$
with $D''v=Dv$ and $D''v'=D'v'$,
$M(F')=(\Lambda (w_1,\cdots ,w_6,v'),D')\otimes (\Lambda u,0)$
and
$M(F'')=(\Lambda (w_1,\cdots ,w_6,v),D)\otimes (\Lambda u,0)$.
%F\times  S^{|w_1|}_{\Q}\times\cdots \times  
%S^{|w_{4}|}_{\Q}$ and  $F''=F\times  S^{|w_5|}_{\Q}\times  
%S^{|w_{6}|}_{\Q}$

(2) 
Consider maps $f:Y\to Z$ and $g:X\to Y$ of rational spaces
whose   homotopy groups are oddly graded.
For even-integers $l,m,n$ with $2\leq l\leq m\leq n$,
put $M(Z)=(\Lambda (w_1,\cdots ,w_n,w),d_Z)$
with $d_Zw_i=0$ and $d_Zw=w_1\cdots w_n$,
 $M(Y)=(\Lambda (w_1,\cdots ,w_n,w,v),D)$
with $Dv=w_1\cdots w_{l}$
and  $M(X)=(\Lambda (w_1,\cdots ,w_n,w,v,u,u'),D')$
with $D'v=w_1\cdots w_{l}$, $D'u=w_1w_2$
and $D'u'=w_l\cdots w_m$.
Then $O(f)=\Q\{w_{l+1},\cdots ,w_n\}$,
 $O(g)=\Q\{w_{3},\cdots ,w_l\}$ and
 $O(f\circ g)=\Q\{w_{3},\cdots ,w_l,w_{m+1},\cdots ,w_n
\}$.
Thus $o(f)=n-l$, 
$o(g)=l-2$, 
$o(f\circ g)=l-m+n-2$
and  in particular $o(f)+o(g)-o(f\circ g)=m-l$ can be arbitrarily  large. 

%(2) Put $M(B)=(\Lambda (w_1,\cdots ,w_{2n}),0)$,
%$Dv=w_1w_2$ and $D'v'=w_1\cdots w_{2n}$ ($n\geq 1$).
%Then $o(f)=o(g)=o(g')=o(g\circ f')=0$ and $o(f')=2n-2$.
%In the above, we can see the examples that
%the `inequality' in Theorem 1.2 is the `equality'. 
%But the author does not know 
%a general condition for $f$ and $g$ to satisfy $o(f\circ g)=o(f)+o(g)$. 
}
\end{exmp}

%\vspace{0.5cm}

\begin{exmp}{\rm
For the homotopy set $[E,B]$ of based maps from $E$ to $B$,
define the subset
$ {\bf G}'(E,B):=\{[f]\in [E,B]\ |\ \ f \mbox{ is a  G-map}\}.$
%$$ GM(X,Y)_{\Q}:=\{[f]\in [X,Y]\ |\ \ f \mbox{ is a  r.G-map}\}$$
A map $f$
from $E$ to $B$ is
said to be a cyclic map if  
$(f|1): E\vee B\to B$ admits an extension
$F:E\times B\to B$ \cite{V}.
The  set of homotopy classes of 
cyclic maps $f:E\to B$ is denoted as ${\bf G}(E,B)$.
Since a cyclic map is a G-map
from $Im\ \pi_*(f)\subset G(B)$ \cite[Lemma 2.1]{Si}(\cite[Corollary 2.2]{LS2}), 
 there is an inclusion 
${\bf G}(E,B)\subset {\bf G}'(E,B)$. 
The quotient map $f :G\to   G/K$ 
for a Lie group $G$ and any closed subgroup $K$ is a cyclic map
 \cite{Si}.
Also the Hopf map $\eta :S^3\to S^2$ is a  cyclic map.
 From \cite[Theorem 2.1]{LMW},
 the map $\eta$ induces $\pi_n(S^3)\cong G_n(S^2)$
 for all $n$.
 Therefore,
if a space $E$ is 2-connected, then any map $f:E\to S^2$
is a G-map.
%For example,$\pi :SO(3)\times S^1\to B=i(S^1)SO(3)\times S^1$
%is a G-map though $B$ is not a for a certain inclusion $S^1\to SO(3)\times S^1$ \cite{Si}. 
 %Notice that it induces $\pi (f):\pi_*(E)\to G(S^2)$.
 %The condition that $\pi (f)(\pi_*(E))\subset G(B)$is  stronger  than it of G-map.
A Gottlieb map is not a  cyclic map in general.
For example, the identity map 
$S^{2n}\overset{=}{\to} S^{2n}$ is not a cyclic map 
\cite[Theorem 3.2]{LS2} but of course a G-map.
%For $n\neq 1,2,3$,$S^{2n+1}\overset{=}{\to} S^{2n+1}$ is not cyclic map but G-map since only $S^3,S^5,S^7$ are H-spaces.???? The collapsing map $CP^n\to S^{2n+1}$
In general, a self-equivalence map $f:B\overset{\simeq}{\to} B$ is a G-map.
  %Let $B$ be an $E$-space for a Lie group $E$.
  %Then the orbit map $f:E\to B$ of the action is a cyclic map.
    We note that a cyclic map factors through an H-space,
 which entails numerous consequences for a cyclic map \cite{LS2}.
But, for our G-map, it seems difficult to search
 such a  useful property. 
%Also see \cite{FL} for cyclic maps.
%Then the model is given by$$(\Lambda W,d_b)\to (\Lambda W\otimes \Lambda (z,V),D)\to (\Lambda (z,V),d)$$ with $H^+(\Lambda W\otimes \Lambda V)=0$.
%Then there is an element $v\in V$ with $Dv=z-w+\alpha$for some $\alpha\in \Lambda W\otimes \Lambda V$with $D\alpha=d_Bw$.Then $\delta_D (w,1)=\sum_i(u_i,f_i)+(v,1)$and $\delta_D((w,1)+\sigma )\neq 0$ for any $\sigma$from $(*)$.

(1)
When 
$H^*(B;\Q)\cong \Q[w]/(w^{k+1})$ with $|w|=2n$,
recall that ${\bf G}(E_{\Q},B_{\Q})\cong H^{2n(k+1)-1}(E;\Q)$ \cite[Example 4.4]{LS2}.
On the other hand, ${\bf G}'(E_{\Q},B_{\Q})\cong [E_{\Q},B_{\Q}]\cong A\times H^{2n(k+1)-1}(E;\Q)$
where $A=\{a\in H^{2n}(E;\Q)| 
a^{k+1}= 0\}$.
%Here G-lement 

%For example,
%consider $[E, B_{\Q}]$
%for $M(E)=\Lambda (v_1,..,v_4,v),d_E)$
%with $|v_i|=3$
%and $H^*(B;\Q )=\Q[w]/(w^2)$ with $|w|=6$.
%Put $M(S^2\times S^5)=(\Lambda (x,y,z),D')$with $|x|=2,|y|=3,|z|=5$, 
%$D'y=x^2$ $D'z=0$ and 
%$(\Lambda W,d_B)=(\Lambda (w,u),d_B) $.
%Then in general, for some $a\in\Q$,
%we can put $$M(f)(w)=av_1..v_4,\ M(f)(u)=av$$

(2) When $H^*(B;\Q )\cong \Q[w]\otimes \Lambda (x,y)/(wxy+w^5)$
with $|w|=2,|x|=3,|y|=5$. 
Then
$B$ is a cohomological symplectic  space %\cite[Def.3.11]{TO}
with formal dimension 16
where $M(B)\cong (\Lambda W,d_B)\cong 
(\Lambda (w,x,y,u),d_B)$ with $|u|=9$, $d_Bw=d_Bx=d_By=0$ and 
$d_Bu=wxy+w^5$.
Put 
$E=S^3\times S^5\times S^9$; i.e., 
$M(E)\cong (\Lambda (v_1,v_2,v_3),0)$ with $|v_1|=3,|v_2|=5,|v_3|=9$.
From  degree arguments we can put
$M(f)(w)=0$, $M(f)(x)=av_1$,  $M(f)(y)=bv_2$, $M(f)(u)=cv_3$
for some $a,b,c\in \Q$. %Then there are the following types (i)$\sim $(v).
Note that, if 
 $a\neq 0,\ b\neq 0$ and $c\neq 0$, it is rational homotopy equivalent
to
the $S^1$-fibration $S^1\to E\to B\simeq ES^1\times_{S^1} E$,
where the model is 
$(\Lambda W,d_B) \to (\Lambda W\otimes \Lambda v,D)\to (\Lambda v,0)$ 
with $|v|=1$ and
 $Dv=w$ \cite{H}.
We see that 
$f$ is an r.G-map if and only if  $a=b=0$ since  
$G(E)_{\Q}=\Q\{v_1^*,v_2^*,v_3^*\}$ and 
$G(B)_{\Q}=\Q\{u^*\}$.
Thus 
$[E_{\Q},B_{\Q}]\cong \Q\times \Q\times \Q$  by $f_{\Q}\equiv (a,b,c)$
and ${\bf G}'(E_{\Q},B_{\Q})=
%(\Q\vee \Q)\times \Q$ and  
{\bf G}(E_{\Q},B_{\Q})\cong \Q$ by $f_{\Q}\equiv (0,0,c)$.}
\end{exmp}

\begin{exmp}{\rm
Put $P_n(Y)$ the $n$th center of the homotopy Lie algebra $\pi_*(\Omega Y)$;
i.e.,
the subgroup of elements $a$ in $\pi_n(Y)$ with $[a,b]=0$ (Whitehead product)
for all $b\in \pi_*(Y)$.
A space $Y$ is called a W-space if $P_n(Y)=\pi_n(Y)$ for all $n$
\cite[Definition 1.8(b)]{Si}. 
\vspace{1mm}

\noindent
{Definition C.}\  {We will call a map $f:E\to B$  a {\it W-map} if 
$\pi_n(f)P_n(E)\subset P_n(B)$ for all $n$.}
\vspace{1mm}

\noindent
For example,
if $\pi_*(f)$ is surjective, $f$ is a W-map.
In spaces, there 
are the implications:
`H-space $\Rightarrow$ {G-space} $\Rightarrow$ { W-space}'
 \cite{Si}.
But `{G-map} $\Rightarrow$ { W-map}'
is false in general.
For example, put
$M(B)=(\Lambda (w_1,w_2,u),d_B)$
with $|w_i|$ odd, $d_Bw_i=0$ and $d_Bu=w_1w_2$.
If the KS-extension $M(B)\to ( \Lambda (w_1,w_2,u,v_1,v_2,v_3,v_4),D)$
of a map $f:E\to B$ is given by $Dv_1=Dv_2=0$, $Dv_3=w_1$ and $Dv_4=w_2uv_1v_2$,
then $f_{\Q}$ is a G-map but not a W-map since 
$w_2^*\not\in G_*(E_{\Q})$ but $w_2^*\in P_*(E_{\Q})$
and $P_*(B_{\Q})=\Q\{u^*\}$.
}
\end{exmp}

%\vspace{0.5cm}

%\vspace{0.3cm}

 \begin{exmp}{\rm  For a fibration, D.Gottlieb proposed a question:
{\it Which homotopy equivalences of the fiber into itself 
  can be extended to fiber homotopy equivalences of the total
  space into itself ?} \cite[\S 5]{G2}.
    We consider a question:
  Which  map $f:E\to B$ can be extended to a map 
  between  fibrations over a sphere,
  that is,
for a fibrations $\xi :E\to E'\to S^{n+1}$,
does there exist
a fibration $\eta :B\to B'\to S^{n+1}$ and 
 a map $f':E'\to B'$ such that the diagram
  $${ (i)\ \ \ \ \ 
\xymatrix{ E\ar[d]_{f}\ar[r]^{}&
 E'\ar@{.>}[d]_{f'}\ar[r]^{}&S^{n+1}\ar@{=}[d]\\
B\ar[r]^{}& 
 B'\ar[r] & S^{n+1}
\\
}}$$
  homotopically commutes ?
% For example,  there is an associated fibration
% ${\mathbb{C}}P^n\to B'\to Y$ satisfying $(i)$  for any odd sphere fibration
%  $S^{2n+1}\to E'\to Y$ and the Hopf map $f:S^{2n+1}\to {\mathbb{C}}P^n$,
%  which is a G-map.
%  Suppose that $E$ and $B$ are CW complexes of finite type.
 If $f:E\to B$ is extended to a map between
 $\xi$ and $\eta$, %in the homotopy exact sequences of them,
 from the result \cite{G2} of Gottlieb, 
 we have a commutative diagram  for all $n$  
    $${ (ii)\ \ \ \ \  
\xymatrix{\pi_{n+1}(S^{n+1})\ar[d]_{\partial_{n+1}^{\ \eta}}
\ar[r]^{\partial_{n+1}^{\ \xi}}&
G_{n}(E)\ar[d]^{\pi_{n} (f)}\\
G_n(B)\ar[r]^{\subset }& 
\pi_{n}(B), 
\\
}}$$
where $\partial_{n+1}$ is the $n+1$th connecting homomorphisms  
in the long exact homotopy sequence of fibration. 
Therefore we have
%To be a  G-map is a necessary condition for a solution to the above problem.

\vspace{1mm}
\noindent
{Claim:}\ 
{\it If $f:E\to B$ is not a G-map,
then there is an $E$-fibration over a sphere where $f$ can not be extended
 to the map
 $f'$ satisfying %any $B$-fibration as 
 $(i)$.}
  \vspace{1mm}

In fact, 
suppose that $0\neq \pi_n(f)(x)\not\in G_n(B)$ for some $x\in G_n(E)$.
Then 
there is a non-trivial fibration $\xi_x :E\to E'\to S^{n+1}$
with $\partial_{n+1}^{\xi_x}(y)=x$ for the generator 
 $y$ of  $\pi_{n+1}(S^{n+1})$  \cite[Thorem I.2]{La}.
Here $\xi_x$ is 
constructed as follows (\cite[page 11]{O}).
Choose 
a preimage $\hat x$ of $x$ under the evaluation map
$\pi_n(aut_1 E)\to G_n(E)$.
From $\pi_{n+1}(Baut_1E)\cong \pi_n(aut_1E)$,
we may consider ${\hat{x}}\in \pi_{n+1}(Baut_1E)$
with representative $S^{n+1}\to Baut_1E$.
Pull back  the universal fibration
over this map to get $\xi_x$.
On the other hand, %even ifthere exists
for any $B$-fibration
$\eta$ over $S^{n+1}$,
 $G_n(B)\ni \partial_{n+1}^{\eta}(y)\neq \pi_n(f)(x)$
 from the assumption.
 Therefore  $(ii)$ does not commute.

But to be a  G-map is not sufficient for the above extension problem.
Let 
$f:E=S^3\times S^5\to S^5=B$
be the projection given by $f(a,b)=b$.
Evidently this is a G-map.
Suppose that a fibration $\xi :E\to E'\to S^3$
is given by a classifying map $h$ with $\pi (h)_{\Q}:
\pi_3(S^3)_{\Q}\cong \pi_3(Baut_1S^3\times S^5)_{\Q}$.
Then the KS-extension of $\xi$ is given by
$(\Lambda (w),0)\to (\Lambda (w,v,v'),D)\to (\Lambda (v,v'),0)$
with $D(v)=0$, $D(v')=wv$,  
$|w|=3$, $|v|=3$ and $|v'|=5$ \cite[Theorem 3.2]{LS}.
Then for any fibration $\eta :B=S^5\to B'\to S^3$,
there is not a map
$f'$ that satisfies $(i)$ 
since $\eta$ is rationally trivial
from degree arguments.
}
\end{exmp}

%\vspace{0.3cm}

\begin{exmp}{\rm
In Example 3.2(1), we see an example of  ``non-Gottlieb'' map
whose homotopy fibre $X$ has the rational homotopy type of  an odd sphere $S^{2n+1}$.
But 
if the homotopy fibre $X$ has the rational homotopy type of  an even  sphere
$S^{2n}$,
then a map $f$ is an r.G-map.
Indeed, 
put $M(S^{2n})=(\Lambda (x,y),d)$
with $|x|=2n$, $|y|=4n-1$, $dx=0$ and $dy=x^2$.
We know that $Dx=0$ and $Dy=x^2+ax+b$ for some 
$a,b\in \Lambda W$
in a KS-extension.
Suppose  $w\in W$ and  $w^*\not\in G(B_{\Q})$.
Then we have, for any $\alpha_i\in \Lambda W$ and 
$\beta_i\in \Lambda W\otimes \Lambda^+(x,y)$ with 
$W=\Q\{w_i\}_{i\in I}$,
$$\delta_f((w,1)+(\sum_{i\in I}w_i,\alpha_i))
\neq \delta_f((\sum_{i\in I}w_i,\beta_i))$$
in $Der_{|w|}(\Lambda W,\Lambda V\otimes \Lambda W)$.
It deduces  $G(B_{\Q},E_{\Q};f_{\Q})\subset G(B_{\Q})$
from Theorem 2.1
and  then $f$ is an r.G-map from Lemma 1.1.

 Recall that
an elliptic space is  one whose rational homology and rational
homotopy are both finite dimensional
and that  an elliptic space $X$ is said to be an $F_0$-space if
the Euler characteristic is positive \cite{FHT}.
When $X$ is an $F_0$-space, for some even degree elements
$x_1,..,x_l$,
there is an isomorphism 
$H^*(X;\Q)\cong \Q[x_1,\cdots ,x_l]/(f_1,\cdots ,f_l)$
with a regular sequence $(f_1,\cdots ,f_l)$
in $ \Q[x_1,\cdots ,x_l]$; i.e.,
$gf_i\in (f_1,\cdots ,f_{i-1})$
implies $g\in (f_1,\cdots ,f_{i-1})$ for any $g\in
\Q [x_1,\cdots ,x_l]$ and all $i$.
For example, $S^{2n}$ is an $F_0$-space with $H^*(S^{2n};\Q)
\cong \Q[x]/(x^2)$.
For an $F_0$-space $X$, S.Halperin conjectures that $Dx_i=0$ for $i=1,..,l$,
which deduces a  fibration with fibre $X$ is totally non-cohomologous to zero \cite{FHT}.
For example, it holds when $X$ is a  homogeneous space \cite{ST}.
%Finally we propose 
%
%\noindent
%{\bf Question A.} 
{
If the homotopy fibre $X$ of a map $f$ is an $F_0$-space,
then is  $f$ an r.G-map ?}
}
\end{exmp}

\vspace{4mm}
\noindent{\bf Acknowledgements.} 
  \ \ It is a pleasure to  thank someone for his kind 
and useful suggestions  for the first version of this papper.
He informed me   the works of \cite{O2}, \cite{O3}, \cite{H2}, \cite{FL}
and indicated that they are 
closely related to the splitting of Theorem 1.3.
Also he pointed out a mistake  in 
Lemma 2.3.

\end{document}